\documentclass[12pt]{amsart}
\usepackage{graphicx}
\usepackage{a4wide}
\usepackage{amssymb}
\usepackage{amsthm}
\usepackage{amsmath}
\usepackage{amscd} \usepackage{verbatim}
\usepackage[all]{xy}
\textheight8.6in \textwidth6.6in \numberwithin{equation}{section}
\theoremstyle{plain}
\newtheorem{theorem}{Theorem}[section]
\newtheorem{corollary}[theorem]{Corollary}
\newtheorem{lemma}[theorem]{Lemma}
\newtheorem{proposition}[theorem]{Proposition}

\theoremstyle{definition}

\newtheorem*{example}{Example}

\theoremstyle{remark}
\newtheorem*{remark}{Remark}
\newtheorem*{notation}{Notation}

\newcommand{\SL}{\text {\rm SL}}
\newcommand{\R}{\mathbb{R}}

\newcommand{\Q}{\mathbb{Q}}

\newcommand{\Z}{\mathbb{Z}}

\newcommand{\C}{\mathbb{C}}
\newcommand{\h}{\mathbb{H}}
\renewcommand{\H}{\mathbb{H}}

\newcommand{\D}{\mathbb{D}}

\newcommand{\zxz}[4]{\begin{pmatrix} #1 & #2 \\ #3 & #4 \end{pmatrix}}

\newcommand{\kzxz}[4]{\left(\begin{smallmatrix} #1 & #2 \\ #3 & #4\end{smallmatrix}\right) }
\newcommand{\kabcd}{\kzxz{a}{b}{c}{d}}

\newcommand{\calF}{\mathcal{F}}

\newcommand{\calM}{\mathcal{M}}

\newcommand{\calO}{\mathcal{O}}

\newcommand{\calQ}{\mathcal{Q}}

\newcommand{\calS}{\mathcal{S}}

\newcommand{\frake}{\mathfrak e}

\newcommand{\Tr}{\text {\rm Tr}}

\newcommand{\bs}{\backslash}

\newcommand{\tr}{\operatorname{tr}}

\newcommand{\Span}{\operatorname{span}}

\newcommand{\Sl}{\operatorname{SL}}
\newcommand{\Gl}{\operatorname{GL}}
\newcommand{\Symp}{\operatorname{Sp}}
\newcommand{\Spin}{\operatorname{Spin}}

\newcommand{\Mp}{\operatorname{Mp}}
\newcommand{\Orth}{\operatorname{O}}

\newcommand{\Mat}{\operatorname{Mat}}

\newcommand{\sig}{\operatorname{sig}}

\newcommand{\SO}{\operatorname{SO}}

\begin{document}

\title[Algebraic formulas for the coefficients of weak  Maass forms] {Algebraic formulas for the coefficients of half-integral weight harmonic weak Maass forms}

\author{Jan Hendrik Bruinier and Ken Ono}
\address{Fachbereich Mathematik, Technische Universit\"at Darmstadt,
Schlossgartenstrasse 7, D-64289, Darmstadt, Germany}
\email{bruinier@mathematik.tu-darmstadt.de}

\address{Department of Math and Computer Science,
Emory University, Atlanta, Georgia 30322}
\email{ono@mathcs.emory.edu}

\thanks{The authors thank the American Institute of Mathematics for their
generous support. The first author thanks the support of DFG grant BR-2163/2-1, and
the second author  thanks the NSF,  the
Hilldale Foundation, the Manasse family, and
the Candler Fund for their support.}

\begin{abstract}
We prove that the coefficients of certain weight $-1/2$ harmonic Maass forms
are ``traces'' of {\it singular moduli} for weak Maass forms. To prove this theorem,
we construct a {\it theta lift} from spaces of weight $-2$ harmonic weak Maass forms to
spaces of weight $-1/2$ vector-valued harmonic weak Maass forms on $\Mp_2(\Z)$, a result
which is of independent interest.
We then prove a general theorem which guarantees (with bounded denominator) when
such Maass  singular moduli are algebraic.
As an example of these  results, we
 derive a formula for the partition function $p(n)$ as a {\it finite} sum
of {\it algebraic} numbers which lie in the usual discriminant $-24n+1$ ring class field.
We indicate how these results extend to general weights. In particular,
we illustrate how one can compute theta lifts for general weights by making use
of the Kudla-Millson kernel and Maass differential operators.
\end{abstract}

\maketitle

%%%%%%%%%%%%%%%%%%%%%%%%%%%%%%
\section{Introduction and statement of results}
%%%%%%%%%%%%%%%%%%%%%%%%%%%%%%

A \textit{partition} \cite{Andrews} of a positive integer $n$ is any nonincreasing sequence
of positive integers which sum to $n$. The partition function $p(n)$, which counts the
number of partitions of $n$,
defines the rapidly increasing
sequence of integers:
\begin{displaymath}
\begin{split}
&1, 1, 2, 3, 5,\dots, p(100)=190569292,\dots,
 p(1000)=24061467864032622473692149727991,\dots.
\end{split}
\end{displaymath}
In celebrated  work \cite{HardyRamanujan}, which gave birth to the ``circle method'', Hardy and Ramanujan
quantified this  rate of growth. They proved the asymptotic formula:
\begin{displaymath}
p(n)\sim \frac{1}{4n\sqrt{3}}\cdot e^{\pi \sqrt{2n/3}}.
\end{displaymath}
Rademacher \cite{Rademacher1, Rademacher2} subsequently perfected this method to derive his famous ``exact'' formula
\begin{equation}\label{exact}
p(n)= 2 \pi (24n-1)^{-\frac{3}{4}} \sum_{k =1}^{\infty} \frac{A_k(n)}{k}\cdot
I_{\frac{3}{2}}\left( \frac{\pi \sqrt{24n-1}}{6k}\right).
\end{equation}
Here $I_{\frac{3}{2}}(\cdot)$ is a  modified Bessel function of the first kind, and $A_k(n)$ is a Kloosterman sum.
%%$$
%%A_k(n):=\frac{1}{2} \sqrt{\frac{k}{12}} \sum_{\substack{x \pmod {24k}\\
%%x^2 \equiv -24n+1 \pmod{24k}}}  (-1)^{\left\{\frac{x}{6}\right\}} \cdot
%%\exp\left(\frac{2\pi i x}{12k}\right).
%%$$
%%In the sum
%%$\{\alpha\}$ denotes the integer nearest to $\alpha$.

\begin{remark} Values of $p(n)$ can be obtained by
rounding sufficiently accurate truncations of (\ref{exact}). Bounding the error between $p(n)$ and such truncations is a well known difficult problem.
Recent work by Folsom and Masri \cite{FolsomMasri} give the best known nontrivial bounds
on this problem.
\end{remark}

We obtain a new formula for $p(n)$.
Answering Questions 1 and 2 of \cite{BringmannOnoPAMS}, we express
$p(n)$ as a {\it finite sum} of {\it algebraic numbers}. These numbers are {\it singular moduli}
for a weak Maass form which we describe using Dedekind's
eta-function $\eta(z):=q^{\frac{1}{24}}\prod_{n=1}^{\infty}(1-q^n)$
(note. $q:=e^{2\pi i z}$ throughout), and the quasimodular Eisenstein series
\begin{equation}
E_2(z):=1-24\sum_{n=1}^{\infty}\sum_{d\mid n}dq^n.
\end{equation}
To this end, we
define the $\Gamma_0(6)$ weight -2 meromorphic modular form $F(z)$ by
\begin{equation}
\label{eq:defF}
F(z):=\frac{1}{2}\cdot \frac{E_2(z)-2E_2(2z)-3E_2(3z)+6E_2(6z)}{\eta(z)^2\eta(2z)^2\eta(3z)^2\eta(6z)^2}
\phantom{:}=q^{-1}-10-29q-\dots.
\end{equation}
Using  the convention that $z=x+iy$, with $x,y\in \R$, we define the {\it weak Maass form}
\begin{equation}
P(z):=-\left(\frac{1}{2\pi i}\cdot \frac{d}{dz}+\frac{1}{2\pi y}\right) F(z)
= \left(1-\frac{1}{2\pi y}\right)q^{-1}+\frac{5}{\pi y}+
\left(29+\frac{29}{2\pi y}\right)q+\dots.
%%\left(208+\frac{52}{\pi \im(z)}\right)q^2+\dots
%%\end{split}
\end{equation}
This nonholomorphic  form has weight 0, and is a weak Maass form
(for more on weak Maass forms, see \cite{BrF1}). It has eigenvalue $-2$ with respect to the hyperbolic Laplacian
\begin{displaymath}
\Delta:=-y^{2}\left(\frac{\partial^2}{\partial x^2}+\frac{\partial^2}{\partial y^2}\right).
\end{displaymath}

To describe our formula, we use discriminant $-24n+1=b^2-4ac$ positive definite integral binary
quadratic forms $Q(x,y)=ax^2+bxy+cy^2$ with the property that
$6\mid a$. The group $\Gamma_0(6)$ acts on such forms, and
we
let $\calQ_n$ be any set of representatives of those equivalence classes with
$a>0$ and $b\equiv 1\pmod{12}$.
For each $Q(x,y)$, we let $\alpha_Q$ be the CM point in $\H$,
the upper half of the complex
plane, for which $Q(\alpha_Q,1)=0$. We then define the ``trace''
\begin{equation}
\Tr(n):=\sum_{Q\in \calQ_n} P(\alpha_Q).
\end{equation}

The following theorem gives the finite algebraic formula for $p(n)$.

\begin{theorem}\label{main}
If $n$ is a positive integer, then
we have that
\begin{displaymath}
   p(n)=\frac{1}{24n-1}\cdot \Tr(n).
\end{displaymath}
The numbers $P(\alpha_Q)$, as $Q$ varies over  $\calQ_n$,
form a multiset of algebraic numbers
which is the union
of Galois orbits for the discriminant $-24n+1$ ring class field.
Moreover, for each $Q\in \calQ_n$ we have that $6(24n-1)P(\alpha_Q)$
is an algebraic integer.
\end{theorem}

Theorem~\ref{main} gives an algorithm for computing $p(n)$, as well as the polynomial
\begin{equation}
H_n(x)=x^{h(-24n+1)}-(24n-1)p(n)x^{h(-24n+1)-1}+\dots:=\prod_{Q\in \calQ_n}(x-P(\alpha_Q))\in \Q[x].
\end{equation}
One simply computes
sufficiently precise approximations of the singular moduli
$P(\alpha_Q)$.

\begin{remark}
Using the theory of Poincar\'e series and identities and formulas for Kloosterman-type sums, one can use Theorem~\ref{main} to give a new (and longer) proof of the exact formula (\ref{exact}).
\end{remark}

\begin{example}
We give an amusing proof of the fact that $p(1)=1$.
In this case, we have that $24n-1=23$, and
we use the $\Gamma_0(6)$-representatives  $$\calQ_1=\{Q_1, Q_2,Q_3\}=\{ 6x^2+xy+y^2, 12x^2+13xy+4y^2, 18x^2+25xy+9y^2\}.$$
The corresponding CM points are
\begin{displaymath}
\alpha_{Q_1}:=-\frac{1}{12}+\frac{1}{12}\cdot \sqrt{-23},\ \ \ \ \
\alpha_{Q_2}:=-\frac{13}{24}+\frac{1}{24}\cdot \sqrt{-23},\ \ \ \ \
\alpha_{Q_3}:=-\frac{25}{36}+\frac{1}{36}\cdot \sqrt{-23}.
\end{displaymath}
Using the explicit Fourier expansion of $P(z)$, we find that $P(\alpha_{Q_3})=\overline{P(\alpha_{Q_2})}$, and that
\begin{displaymath}
\begin{split}
P(\alpha_{Q_1})\sim 13.965486281 \ \ \ \ \ {\text {\rm and}}\ \ \  \ \
P(\alpha_{Q_2})\sim  4.517256859-3.097890591 i.
\end{split}
\end{displaymath}
By means of these numerics we can prove that
\begin{displaymath}
H_1(x):=\prod_{m=1}^{3}(x-P(\alpha_{Q_m}))=x^3-23x^2+\frac{3592}{23}x-419,
\end{displaymath}
and this confirms that $p(1)=\frac{1}{23}\Tr(1)=1$.
If $\beta:=161529092+18648492\sqrt{69}$, then we have
\begin{displaymath}
\begin{split}
P(\alpha_{Q_1})&=\frac{\beta^{1/3}}{138}+\frac{2782}{3\beta^{1/3}}+\frac{23}{3},\\
P(\alpha_{Q_2})&=-\frac{\beta^{1/3}}{276}-\frac{1391}{3\beta^{1/3}}+\frac{23}{3}-\frac{\sqrt{-3}}{2}\cdot
\left(\frac{\beta^{1/3}}{138}-\frac{2782}{3\beta^{1/3}}\right).
%%P(\alpha_{Q_3})&=-\frac{\beta^{1/3}}{276}-\frac{1391}{3\beta^{1/3}}+\frac{23}{3}+\frac{\sqrt{-3}}{2}\cdot
%%\left(\frac{\beta^{1/3}}{138}-\frac{2782}{3\beta^{1/3}}\right).
\end{split}
\end{displaymath}
\end{example}

The claim in Theorem~\ref{main} that $p(n)=\Tr(n)/(24n-1)$ is an example of a general theorem (see Theorem~\ref{GeneralThm})
on  ``traces'' of CM values of certain weak Maass forms. This result pertains to
 weight $0$ weak Maass forms which are images under the Maass raising operator of
weight $-2$ {\it harmonic Maass forms}. We apply this to $F(z)$
which is a weight $-2$ {\it weakly holomorphic modular form}, a meromorphic modular form
whose poles are supported at cusps. Theorem~\ref{GeneralThm}
is a new result which adds to the extensive literature
(for example, see \cite{BringmannOno, BrF,
BrJO, DIT1, DIT2, DukeJenkins, Jenkins, millerpixton, cbms}) inspired by  Zagier's seminal paper \cite{Zagier} on ``traces'' of singular moduli.

To obtain this result, we employ the theory of theta lifts as in earlier work by Funke and the first author
\cite{BrF1, Funke}. Here we use the Kudla-Millson theta functions  to construct a new theta
lift (see Corollary~\ref{cor:lift}),  a result which is of independent interest. The lift maps spaces of weight $-2$ harmonic weak Maass forms
to spaces of weight $-1/2$ vector valued harmonic weak Maass forms for $\Mp_2(\Z)$. In Section~\ref{theta} we recall properties of
these theta functions, and in Section~\ref{General} we construct the lift, and we then employ an argument of Katok and Sarnak to prove
Theorem~\ref{GeneralThm}. In Section~\ref{General} we also indicate how Corollary~\ref{cor:lift} and Theorem~\ref{GeneralThm} extend to
general weights. In particular, we illustrate how to define theta lifts for general
weights using the Kudla-Millson kernel and Maass differential operators. For the sake
of brevity and the application to $p(n)$, we chose to focus on the case of weight $-2$
harmonic Maass forms.

To complete the proof of Theorem~\ref{main}, we must show that the values in $P(\alpha_Q)$ are  {\it algebraic numbers} with {\it bounded denominators}.
To prove these claims, we require the classical theory of
complex multiplication,  as well as new results which bound denominators
of suitable singular moduli. For example, we bound the denominators of
the singular moduli
(see Lemma~\ref{lem:bound2})
 of suitable nonholomorphic modular functions which contain the nonholomorphic Eisenstein series $E_2^*(z)$ as a factor. Theorem~\ref{thm:boundgen} is our general result
 which bounds the denominators of algebraic  Maass singular moduli such as $P(\alpha_Q)$. These results are contained in Section~\ref{CM}.
In Section~\ref{examples} we give
further examples, and we then conclude with some natural questions.

%%%%%%%%%
%%\section*{Acknowledgements} %%%%%%%%%
%%%%%%\putin

%%%%%%%%%%%%%%%%%%%%%%%%%%%%%%
\section{The Kudla-Millson theta functions}\label{theta}
%%%%%%%%%%%%%%%%%%%%%%%%%%%%%%
We assume that the reader is familiar with basic facts about weak Maass forms
(for example, see \cite{BrF}).
Using the Kudla-Millson theta functions, we will construct a theta lift from spaces of weight $-2$ harmonic weak Maass forms on
$\Gamma_0(N)$
to weight $-1/2$ vector-valued harmonic Maass forms on $\Mp_2(\Z)$.
This lift will be crucial to the proof of Theorem~\ref{GeneralThm} which interprets coefficients of holomorphic parts of these
weight $-1/2$ forms as ``traces'' of the CM values of weight 0 weak Maass forms
with Laplacian eigenvalue $-2$.

We begin by recalling some important facts
 about these theta functions in the setting of the present paper (see \cite{KMI}, \cite{BrF}).
Let $N$ be a positive integer.
Let $(V,Q)$ be the quadratic space over $\Q$ of signature $(1,2)$
given by the trace zero $2\times 2 $ matrices
\begin{equation} \label{iso}
V  :=\left\{ X = \begin{pmatrix} x_1 & x_2 \\ x_3 & -x_1
 \end{pmatrix} \in \Mat_2(\Q) \right\},
\end{equation}
with the quadratic form $Q(X) = N\det(X)$. The corresponding bilinear
form is $(X,Y) = -N\tr(XY)$.  We let $G = \Spin(V)$, viewed
as an algebraic group over $\Q$, and write $\bar{G}$ for
its image in $\SO(V)$. We realize the associated
symmetric space $\D$ as
the Grassmannian of lines in $V(\R)$ on which the quadratic form $Q$ is
positive definite:
\[
\D \simeq \{ z \subset V(\R) ; \;\text{$\dim z =1$ and $Q|_z > 0$} \}.
\]
The group $\Sl_2(\Q)$ acts on $V$ by conjugation
\[
g.X := gXg^{-1}
\]
for $X\in V$ and $g\in \Sl_2(\Q)$. This gives rise to the isomorphisms
$G\simeq\SL_2$ and $\bar G\simeq\operatorname{PSL_2}$ .

The hermitian symmetric space $\D$ can be identified with the complex upper half plane $\H$
as follows: We choose as a base point $z_0\in {\D}$ the line spanned by
$\left( \begin{smallmatrix} 0 & 1 \\ -1 & 0 \end{smallmatrix}\right)$. Its stabilizer in ${G}(\R)$ is equal to  ${K} = \SO(2)$.
For $z= x+iy
\in \H$, we choose $g_z \in {G}(\R)$ such that $g_zi = z$. If we associate to $z$ the vector
\begin{equation}
X(z) :=  g_z.\zxz{0}{1}{-1}{0} = \frac1{y}
\begin{pmatrix} -x & z\bar{z} \\ -1 & x
\end{pmatrix}\in V(\R),
\end{equation}
then $Q(X(z)) = N$, and $g.X(z) = X(gz)$ for $g \in G(\R)$.
We obtain the isomorphism
\begin{equation}
\H  \longrightarrow\D,\quad
z \mapsto g_z z_0
%= \Span\left(g_z . \left( \begin{matrix} 0 & 1 \\ -1 & 0\end{matrix}\right)\right)
= \R X(z).
\end{equation}
The minimal majorant of $(\;,\;)$
associated to $z \in \D$ is given by $(X,X)_z = (X,X(z))^2 - (X,X)$.

Let $L \subset V$ be an even lattice and write $L'$ for the dual lattice. Let $\Gamma$ be a congruence subgroup of $\Spin(L)$, which takes $L$ to itself and acts trivially on the discriminant group $L'/L$.
We write $M = \Gamma \backslash \D$ for the associated modular curve.

\begin{comment}
\begin{notation}
From now on, we will write $z= x+iy$ for an
element in the orthogonal symmetric space $\D \simeq \h $. The upper case letter $X$  we reserve for vectors in $V(\R)$, thought of as elements in $B_0(\R)$. Its coefficients we denote by $x_i$. Later on, we will write $\tau = u+iv \in \h$ for a modular form variable in $\h$; i.e, we consider $\tau$ as a variable for the (symplectic) symmetric space associated to $\SL_2 \simeq\Symp(1)$.
\end{notation}
\end{comment}

Heegner points in $M$ are given as follows.
If  $X \in V(\Q)$ with $Q(X) >0$, we put
\begin{equation}
\D_X = \Span(X) \in \D.
\end{equation}
The stabilizer $\Gamma_X\subset \Gamma$ of $X$  is finite. We denote by $Z(X)$ the image of $\D_X$ in $M$, counted with multiplicity $\tfrac1{|\bar{\Gamma}_X|}$. 
%Moreover, we put $Z_X = \emptyset$ if $Q(X) \leq 0$.
For $m \in \Q_{>0}$ and $h \in L'/L$,  the group $\Gamma$ acts on
${L}_{m,h} = \{X \in L +h ;\; Q(X) =m\}$
with finitely many orbits. We define the \emph{Heegner divisor} of discriminant $(m,h)$  on $M$ by
\begin{equation}
Z(m,h) = \sum_{X \in \Gamma \bs L_{m,h} }Z(X).
\end{equation}

\subsection{The Kudla-Millson function}

Kudla and Millson defined \cite{KMI} a Schwartz function $\varphi_{KM}$ on $V(\R)$ valued in $\Omega^{1,1}(\D)$, the differential forms on $\D$ of Hodge type $(1,1)$, by
\begin{equation}
\varphi_{KM}(X,z) = \biggl( (X,X(z))^2  - \frac1{2\pi} \biggr) \, e^{-\pi (X,X)_z}
\, \Omega,
\end{equation}
where $ \Omega =   \tfrac{dx \wedge dy}{y^2} = \tfrac{i}2 \tfrac{dz \wedge d\bar{z}}{y^2}$.
We have $\varphi (g.X,g z) = \varphi (X,z)$ for $g \in G(\R)$.
% i.e., the pullback of $\varphi$ with respect to the action of $g$ satisfies $g^{\ast} \varphi(X,z) = \varphi(g^{-1}X,z)$.
We define
\begin{align}
\varphi_{KM}^0(X,z) = e^{\pi(X,X)} \varphi_{KM} (X,z) =\left((X,X(z))^2  - \frac1{2\pi} \right) \, e^{-2 \pi R(X,z) } \, \Omega,
\end{align}
where, following \cite{KAnn}, we set
\begin{equation}\label{Rformel}
R(X,z) := \frac12(X,X)_z - \frac12(X,X) = \frac{1}{2N}(X,X(z))^2 - (X,X).
\end{equation}
The quantity $R(X,z)$ is always non-negative. It equals $0$ if and only if $z =\D_X$, that is,  if $X$ lies in the line generated by $X(z)$. Hence, for $X \ne 0$, this does not occur if $Q(X) \leq 0$.
Recall that for $Q(X) >0$, the $2$-form $\varphi_{KM}^0(X,z)$ is a Poincar\'e dual form for the Heegner point $\D_X$,  while it is exact for $ Q(X) <0$.

\subsection{The Weil representation}

We write $\Mp_2(\R)$ for the metaplectic two-fold cover of
$\Sl_2(\R)$. The elements of this group are pairs $(M,\phi(\tau))$,
where $M=\kabcd\in\Sl_2(\R)$ and $\phi:\H\to \C$ is a holomorphic
function with $\phi(\tau)^2=c\tau+d$.  The multiplication is defined
by
\[
(M,\phi(\tau)) (M',\phi'(\tau))=(M M',\phi(M'\tau)\phi'(\tau)).
\]
We denote the integral metaplectic group, the inverse image of
$\Sl_2(\Z)$ under the covering map, by
$\tilde\Gamma=\Mp_2(\Z)$. It is well known that $\tilde\Gamma$ is
generated by $T:= \left( \kzxz{1}{1}{0}{1}, 1\right)$, and $S:=
\left( \kzxz{0}{-1}{1}{0}, \sqrt{\tau}\right)$.
%One has the
%relations $S^2=(ST)^3=Z$, where $Z:=\left( \kzxz{-1}{0}{0}{-1},
%i\right)$ is the standard generator of the center of $\tilde\Gamma$.
We let $\tilde\Gamma_\infty:=\langle T\rangle\subset\tilde\Gamma$.

We denote the
standard basis elements of the group ring $\C[L'/L]$ by
$\frake_h$ for $h\in L'/L$.
%, and
%write $\langle\cdot,\cdot \rangle$ for the standard scalar product
%(antilinear in the second entry) such that $\langle
%\frake_h,\frake_{h'}\rangle =\delta_{h,h'}$.
Recall (for example, see \cite{Bo1}, \cite{Br}) that the Weil representation
$\rho_L$ associated with the discriminant form $L'/L$ is the unitary
representation of $\tilde\Gamma$ on
%the group algebra
$\C[L'/L]$ defined by
\begin{align}
\label{eq:weilt}
\rho_L(T)(\frake_h)&:=e(h^2/2)\frake_h,\\
\label{eq:weils} \rho_L(S)(\frake_h)&:=
\frac{e(-\sig(L'/L)/8)}{\sqrt{|L'/L|}} \sum_{h'\in L'/L} e(-(h,h'))
 \frake_{h'}.
\end{align}
%\intertext{Note that}
%\begin{align}
%\label{eq:weilz} \rho_L(Z)(\frake_h)&=e((b^--b^+)/4) \frake_{-h}.
%\end{align}
Here $\sig(L'/L)$ denotes the signature of the discriminant form $L'/L$ modulo $8$.

For $k\in \frac{1}{2}\Z$, we let $H_{k,\rho_L}$ be the space of $\C[L'/L]$-valued harmonic Maass forms of weight $k$ for the group $\tilde\Gamma$ and the representation $\rho_L$. We write $M^!_{k,\rho_L}$ for the subspace of weakly holomorphic forms (see \cite{BrF} for definitions). We note that
$H_{k,\rho_L}$ is denoted $H_{k,L}^{+}$ in \cite{BrF}.

\subsection{Theta series}

For $\tau = u+ iv \in \h$ with $u,v\in \R$, we put $g'_{\tau} =
\left( \begin{smallmatrix}1&u\\0&1\end{smallmatrix} \right)
\left( \begin{smallmatrix}v^{1/2}&0\\0&v^{-1/2}\end{smallmatrix}
\right)$.  We denote by $\omega$ the Weil representation of $\Mp_2(\R)$ on  the space of Schwartz functions $S(V(\R))$. For $h\in L'/L$ and $\varphi\in S(V(\R))$
of weight $k$ with respect to the action of $\widetilde{\SO}(2,\R)\subset
\Mp_2(\R)$, we define the $\C[L'/L]$-valued theta function
\begin{align}
%\theta_{h}(\tau,\varphi) &= \sum_{X \in
%h + L} v^{-k/2} (\omega( g'_{\tau}) \varphi)(X) ,\\
\Theta_L(\tau,\varphi)&= \sum_{X\in L'}  v^{-k/2} (\omega( g'_{\tau}) \varphi)(X) \frake_X.
\end{align}
It is well known that
$\Theta_L(\tau,\varphi)$ is a (in general non-holomorphic) modular form of weight $k$ for $\tilde\Gamma$ with representation $\rho_L$.
In particular, for the Kudla-Millson Schwartz function, we obtain, in the variable $\tau$, that
\[
\Theta_L(\tau,z,\varphi_{KM}):= \Theta_L(\tau,z,\varphi_{KM}(\cdot,z))
\]
is a non-holomorphic modular form of weight $3/2$ for $\tilde\Gamma$ with representation $\rho_L$. In $z$ it is a $\Gamma$-invariant $(1,1)$-form on $\D$.

We will also be interested in the standard Siegel theta function
\[
\Theta_L(\tau,z,\varphi_{S}):= \Theta_L(\tau,z,\varphi_{S}(\cdot,z)),
\]
where $\varphi_S(X,z)= e^{-\pi(X,X)_z}$ is the Gaussian on $V(\R)$ associated to the majorant $(\cdot,\cdot)_z$. In $\tau$, it is a non-holomorphic modular form of weight $-1/2$ for $\tilde\Gamma$ with representation $\rho_L$, while it is a $\Gamma$-invariant function in $z$.
Explicitly we have
\[
\Theta_L(\tau,z,\varphi_{S}):= v\sum_{X\in L'} e^{-2\pi vR(X,z)} e( Q(X)\tau )\frake_X.
\]

\subsection{Differential operators}

Let $k\in \frac{1}{2}\Z$.
Recall that the hyperbolic Laplace operator of weight $k$  on functions
in the variable $\tau$ on $\h$ is given by
\begin{equation}\label{defdelta}
\Delta_k=\Delta_{k,\tau} = -v^2\left( \frac{\partial^2}{\partial u^2}+ \frac{\partial^2}{\partial v^2}\right) + ikv\left( \frac{\partial}{\partial u}+i \frac{\partial}{\partial v}\right).
\end{equation}
The Maass raising and lowering operators on non-holomorphic modular forms of weight $k$
are defined as the differential operators
\begin{align}
R_k  &=2i\frac{\partial}{\partial\tau} + k v^{-1} , \\
L_k  &= -2i v^2 \frac{\partial}{\partial\bar{\tau}}.
\end{align}
The raising operator $R_k$ raises the weight of an automorphic form by
$2$, while $L_k$ lowers it by $2$.  The Laplacian $\Delta_k$ can be
expressed in terms of $R_k$ and $L_k$ by
\begin{equation}\label{deltalr}
-\Delta_k = L_{k+2} R_k +k = R_{k-2} L_{k}.
\end{equation}

We let $\partial$, $\bar{\partial}$ and $d$ be the usual differentials
on $\D$. We set $d^c = \tfrac1{4\pi i} (\partial - \bar{\partial})$,
so that $dd^c = - \tfrac1{2\pi i} \partial \bar{\partial}$.  According
to \cite[Theorem 4.4]{BrF1}, the Kudla-Millson theta function and the
Siegel theta function are related by the identity
\begin{align}
\label{eq:diffid}
L_{3/2,\tau} \Theta_L(\tau,z,\varphi_{KM}) = -dd^c \Theta_L(\tau,z,\varphi_{S})
= \frac{1}{4\pi}\Delta_{0,z} \Theta_L(\tau,z,\varphi_{S})\cdot  \Omega.
\end{align}
Moreover, by \cite[Prop. 4.5]{Br}, it follows that the Laplace operators on the Kudla-Millson theta kernel are related by
\begin{align}
\label{eq:laplace}
\Delta_{3/2,\tau} \Theta_L(\tau,z,\varphi_{KM})
= \frac{1}{4}\Delta_{0,z} \Theta_L(\tau,z,\varphi_{KM}).
\end{align}

\subsection{A lattice related to $\Gamma_0(N)$}

For the rest of this section, we let $L$ be the even lattice
\begin{align}
\label{latticeN}
 L:=\left\{\zxz{b}{a/N}{c}{-b}:\quad a,b,c\in
\Z\right\}.
\end{align}
The dual lattice is given by
\begin{align}
\label{latticeN2}
 L':=\left\{\zxz{b/2N}{a/N}{c}{-b/2N}:\quad
\text{$a,b,c\in \Z$} \right\}.
\end{align}
We identify $L'/L$ with $\Z/2N\Z$, and the quadratic form on $L'/L$
is identified with the quadratic form $x\mapsto -x^2$ on $\Z/2N\Z$.
The level of $L$ is $4N$.
The group $\Gamma_0(N)$ is contained in $\Spin(L)$ and acts trivially on $L'/L$.
We denote by $\ell$, $\ell'$ the primitive isotropic vectors
\[
\ell=\zxz{0}{1/N}{0}{0},\qquad \ell'=\zxz{0}{0}{-1}{0}
\]
in $L$, and write $K$ for the one-dimensional lattice $\Z\kzxz{1}{0}{0}{-1}\subset L$.
We have $L=K+\Z\ell+\Z\ell'$ and $L'/L\cong K'/K$.
For $\lambda\in V(\R)$ and $z=x+iy\in \H$ we let $\lambda_z\in V(\R)$ be the orthogonal projection of $\lambda$ to $\R X(z)$. It is easily checked that
\begin{displaymath}
\ell_z= \frac{1}{2Ny} X(z)\ \ \ {\text {\rm and}}\ \ \
\ell_z^2= \frac{1}{2Ny^2}.
%\\
%\frac{\ell_z}{\ell_z^2}&=\zxz{-x}{z\bar z}{-1}{x}.
\end{displaymath}

Following \cite{Bo1} we define a theta function for the smaller lattice $K$  as follows. For $\alpha,\beta\in  K\otimes\R$ and $h\in K'/K$ we put
\begin{align*}
\xi_h(\tau,\alpha,\beta)&= \sqrt{v} \sum_{\lambda\in K+h} e\big( Q(\lambda+\beta)\bar \tau-(\lambda+\beta/2,\alpha)\big),\\
\Xi_K(\tau,\alpha,\beta)&= \sum_{h\in K'/K}\xi_h(\tau,\alpha,\beta)\frake_h.
\end{align*}
According to \cite[Theorem 4.1]{Bo1}, the function $\Xi_K(\tau,\alpha,\beta)$ transforms like a non-holomorphic modular form of weight $-1/2$ for $\tilde\Gamma$ with representation $\rho_K$.
For $z\in \H$ we put $\mu(z)=\kzxz{-x}{0}{0}{x}\in K\otimes\R$.
Theorem 5.2 of \cite{Bo1} allows us to rewrite $\Theta_L$ as a Poincar\'e series:

\begin{proposition}
\label{prop:poincare}
We have that
\[
\Theta_L(\tau,z,\varphi_S) = \frac{1}{\sqrt{2\ell_z^2}} \cdot \Xi_K(\tau,0,0)
+ \frac{1}{2\sqrt{2\ell_z^2}}\sum_{n=1}^\infty \sum_{\gamma\in \tilde\Gamma_\infty\bs \tilde \Gamma}
\left[ \exp\left(-\frac{\pi n^2}{2v\ell_z^2}\right) \Xi(\tau,n\mu(z) ,0)\right]\mid_{-1/2,\rho_K}\gamma.
\]
\end{proposition}

\subsection{Poincar\'e series}

We now recall some facts on Poincar\'e series with exponential growth at the cusps.
Let $k\in \frac{1}{2}\Z$.
Let $M_{\nu,\,\mu}(z)$ and $W_{\nu,\,\mu}(z)$ \label{bi4} be the usual Whittaker functions (see p.~190 of \cite{AS}).
%% \cite{B1} Vol.~I Chap.~6 p.~264.
For convenience, we put for $s\in\C$ and $y\in\R_{>0}$:
\begin{align}\label{calM}
\calM_{s,k}(y)= y^{-k/2} M_{-k/2,\,s-1/2}(y).
\end{align}
For $s=k/2$, we have the identity
\begin{align}\label{Mspecial}
\calM_{k/2,k}(y)&= y^{-k/2}  M_{-k/2,\,k/2-1/2}(y)=e^{y/2}.
\end{align}

Let $\Gamma_\infty$ be the subgroup of $\Gamma=\Gamma_0(N)$ generated by $\kzxz{1}{1}{0}{1}$.
If $k$ is integral, and $m$ is a positive integer, we define the Poincar\'e series
\begin{equation}\label{DefF}
F_{m}(z,s,k)=\frac{1}{2\Gamma(2s)}\sum_{\gamma\in \Gamma_\infty\bs \Gamma} \left[ \calM_{s,k}(4\pi m y) e(-m x)\right]\mid_k\gamma,
\end{equation}
where $z=x+iy\in\H$ and $s\in\C$ with $\Re(s)>1$ (for example, see~\cite{Br} ).
This Poincar\'e series converges for $\Re(s)>1$, and it is an eigenfunction of $\Delta_k$ with eigenvalue
$s(1-s)+(k^2-2k)/4$. Its specialization at $s_0=1-k/2$ is a harmonic Maass form \cite[Proposition 1.10]{Br}.
Its principal part at the cusp $\infty$ is given by $q^{-m}+C$ for some constant $C\in \C$, and the principal parts at the other cusps are constant.

The next proposition describes the images of these series under the Maass raising operator.
\begin{proposition}
\label{prop:raisef}
We have that
\[
\frac{1}{4\pi m} R_k F_{m}(z,s,k) = (s+k/2) F_{m}(z,s,k+2).
\]
\end{proposition}

\begin{proof}
Since $R_k$ commutes with the slash operator, it suffices to show that
\[
\frac{1}{4\pi m} R_k\calM_{s,k}(4\pi m y) e(-m x) = (s+k/2) \calM_{s,k+2}(4\pi m y) e(-m x).
\]
This identity follows from (13.4.10) and (13.1.32) in \cite{AS}.
\end{proof}

We also define $\C[L'/L]$-valued analogues of these series.
Let $h\in L'/L$,  and let $m\in \Z-Q(h)$ be positive.
For $k\in \Z-\frac{1}{2}$ we define
\begin{equation}\label{DefcalF}
\calF_{m,h}(\tau,s,k)=\frac{1}{2\Gamma(2s)}\sum_{\gamma\in \tilde\Gamma_\infty\bs \tilde\Gamma} \left[ \calM_{s,k}(4\pi m y) e(-m x)\frake_h\right]\mid_{k,\rho_L}\gamma.
\end{equation}
The series converges for $\Re(s)>1$ and defines a weak Maass form of weight $k$ for $\tilde \Gamma$ with representation $\rho_L$. The special value at $s=1-k/2$ is harmonic. If $k\in\Z-\frac{1}{2}$, it has the principal part
$ q^{-m}\frake_h+ q^{-m}\frake_{-h}+C$ for some constant $C\in\C[L'/L]$.

\section{The theta lift and ``traces'' of CM values of  weak Maass forms}\label{General}

Here we construct the theta lift which we then use to prove that
the coefficients of certain weight $-1/2$ harmonic weak Maass forms
are ``traces'' of CM values of weak Maass forms.

\subsection{A theta lift}
\label{sect:3.1}

Let $L$ be the lattice \eqref{latticeN}.
For $k\in \frac{1}{2}\Z$, we let $H_k( N)$ denote the space of harmonic Maass forms of
weight $k$ for $\Gamma:=\Gamma_0(N)$. We let $H^\infty_{k}(N)$ denote the subspace of $H_{k}(N)$ consisting of those harmonic Maass forms whose principal parts at all cusps other than $\infty$ are constant.
We write $M_k^!(N)$ for the subspace of   weakly holomorphic forms in $H_{k}(N)$, and we put
$M_k^{!,\infty}(N)=M_k^!(N)\cap H^\infty_{k}(N)$.

%%In this subsection we define a theta lift from
%%weak Maass forms of weight $-2$ for the group $\Gamma$ to weak Maass forms of weight $-1/2$ for $ \tilde \Gamma$ with represesentation $\rho_L$.
%
For a weak Maass form $f$ of weight $-2$ for $\Gamma$ we define
\begin{align}
\Lambda(\tau,f)=  L_{3/2,\tau}\int_{M} (R_{-2,z}f(z)) \Theta_L(\tau,z,\varphi_{KM}).
\end{align}
According to \cite[Proposition 4.1]{BrF}, the Kudla-Millson theta kernel has exponential decay as $O(e^{-Cy^2})$ for $y\to \infty$ at all cusps of $\Gamma$ with some constant $C>0$. Therefore the theta integral converges absolutely. It defines a $\C[L'/L]$-valued function on $\H$ that transforms like a non-holomorphic modular forms of weight $-1/2$ for $\tilde\Gamma$.
We denote by $\Lambda_h(\tau,f)$ the components of the lift $\Lambda(\tau,f)$ with respect to the standard basis $(\frake_h)_h$ of $\C[L'/L]$.

The group $\Orth(L'/L)$ can be identified with the group generated by the Atkin-Lehner involutions. It acts on weak Maass forms for $\Gamma$ by the Petersson slash operator. It also acts on $\C[L'/L]$-valued modular forms with respect to the Weil representation $\rho_L$ through the natural action on $\C[L'/L]$. The following proposition, which is easily checked, shows that the theta lift is equivariant with respect to the action of $\Orth(L'/L)$.

\begin{proposition}
\label{prop:acto}
For $\gamma\in \Orth(L'/L)$ and $h\in L'/L$, we have
\[
\Lambda_{\gamma h}(\tau,f)= \Lambda_{h}(\tau,f\mid_{-2}\gamma^{-1}).
\]
\end{proposition}

\begin{proposition}
\label{prop:eigval}
If $f$ is an eigenform of the Laplacian $\Delta_{-2,z}$ with eigenvalue $\lambda$, then
$\Lambda(\tau,f)$ is an eigenform of $\Delta_{-1/2,\tau}$ with eigenvalue $\lambda/4$.
\end{proposition}

\begin{proof}
The result follows from \eqref{eq:laplace} and the fact that
\begin{align}
\label{eq:comm}
R_{k}\Delta_{k}=(\Delta_{k+2}-k)R_k,\qquad \Delta_{k-2}L_k= L_k(\Delta_k+2-k).
\end{align}
We may use  symmetry of the Laplacian on the functions in the integral because of the very rapid decay of the Kudla-Millson theta kernel \cite[Proposition 4.1]{BrF}.
\end{proof}

We now compute the lift of the Poincar\'e series.

\begin{theorem}
\label{thm:thetamain}
If $m$ is a positive integer, then we have
\[
\Lambda\big(\tau,F_{m}(z,s,-2)\big)=
\frac{2^{2-s} \sqrt{\pi}Ns(1-s)}{\Gamma(\frac{s}{2}-\frac{1}{2})}
\sum_{n\mid m} n\cdot
\calF_{\frac{m^2}{4Nn^2},\frac{m}{n}}(\tau,\tfrac{s}{2}+\tfrac{1}{4},-\tfrac{1}{2}).
\]
\end{theorem}

\begin{proof}
By definition we have
\begin{align*}
\Lambda\big(\tau,F_{m}(z,s,-2)\big)= L_{3/2,\tau}\int_{M} (R_{-2,z}F_{m}(z,s,-2)) \Theta_L(\tau,z,\varphi_{KM}).
\end{align*}
Employing Propsition \ref{prop:raisef} and \eqref{eq:diffid}, we see that this is equal to
\begin{align*}
m(s-1)\int_{M} F_{m}(z,s,0) \Delta_{0,z}\Theta_L(\tau,z,\varphi_{S})\Omega.
\end{align*}
Using definition \eqref{DefF}, we find, by the usual unfolding argument, that
\begin{align*}
\Lambda\big(\tau,F_{m}(z,s,-2)\big)=\frac{m(s-1)}{\Gamma(2s)} \int_{\Gamma_{\infty}\bs \H} \calM_{s,0}(4\pi m y) e(-m x)\Delta_{0,z}\Theta_L(\tau,z,\varphi_{S})\Omega.
\end{align*}
By Proposition \ref{prop:poincare}, we may replace
$\Delta_{0,z}\Theta_L(\tau,z,\varphi_{S})$ by
$\Delta_{0,z} \tilde \Theta_L(\tau,z,\varphi_{S})$, where
\[
\tilde \Theta_L(\tau,z,\varphi_{S})=\frac{1}{2\sqrt{2\ell_z^2}}\sum_{n=1}^\infty \sum_{\gamma\in \tilde\Gamma_\infty\bs \tilde \Gamma}
\left[ \exp\left(-\frac{\pi n^2}{2v\ell_z^2}\right) \Xi(\tau,n\mu(z) ,0)\right]\mid_{-1/2,\rho_K}\gamma.
\]
Recall that $\ell_z^2=\frac{1}{2Ny^2}$. The function $\tilde \Theta_L(\tau,z,\varphi_{S})$ and its partial derivatives
have square exponential decay as $y \to\infty$.
Therefore, for $\Re(s)$ large, we may move the Laplace operator to the Poincar\'e series and obtain
\begin{align}
\label{eq:lifti}
\Lambda\big(\tau,F_{m}(z,s,-2)\big)&=\frac{m(s-1)}{\Gamma(2s)} \int_{\Gamma_{\infty}\bs \H} \big(\Delta_{0,z}\calM_{s,0}(4\pi m y) e(-m x)\big) \tilde\Theta_L(\tau,z,\varphi_{S})\Omega\\
\nonumber
&=-\frac{m s(s-1)^2}{\Gamma(2s)} \int_{\Gamma_{\infty}\bs \H} \calM_{s,0}(4\pi m y) e(-m x) \tilde\Theta_L(\tau,z,\varphi_{S})\Omega\\
\nonumber
&=-\frac{m s(s-1)^2}{\Gamma(2s)} \sum_{n=1}^\infty\sum_{\gamma\in \tilde\Gamma_\infty\bs \tilde \Gamma}  I(\tau,s,m,n)\mid_{-1/2,\rho_K}\gamma,
\end{align}
where
\begin{align*}
I(\tau,s,m,n)=  \int_{y=0}^\infty\int_{x=0}^1 \calM_{s,0}(4\pi m y) e(-m x)
\frac{1}{2\sqrt{2\ell_z^2}}
\exp\left(-\frac{\pi n^2}{2v\ell_z^2}\right) \Xi(\tau,n\mu(z) ,0)\,\frac{dx\,dy}{y^2}.
\end{align*}
If we use the fact that
$K'=\Z\kzxz{1/2N}{0}{0}{-1/2N}$, and identify $K'/K\cong \Z/2N \Z$, then we have
\begin{align*}
 \Xi(\tau,n\mu(z) ,0)=\sqrt{v}\sum_{b\in \Z} e\left(- \frac{b^2}{4N}\bar\tau-nbx\right)\frake_b.
\end{align*}
Inserting this in the formula for $I(\tau,s,m,n)$, and by integrating over $x$, we see that
$I(\tau,s,m,n)$ vanishes when  $n\nmid m$. If $n\mid m$, then  only the summand for $b=-m/n$ occurs, and so
\begin{align*}
I(\tau,s,m,n)&=  \frac{\sqrt{Nv}}{2} \int_{0}^\infty \calM_{s,0}(4\pi m y) \exp\left(-\frac{\pi N n^2 y^2}{v}\right)\frac{dy}{y}
 e\left(- \frac{m^2}{4Nn^2}\bar\tau\right)\frake_{-m/n}.
\end{align*}
To compute this last integral, we note that
\begin{align*}
\calM_{s,0}(4\pi m y)= M_{0,s-1/2}(4\pi m y) =  2^{2s-1}\Gamma(s+1/2)\sqrt{4\pi m y}\cdot I_{s-1/2}(2\pi m y)
\end{align*}
(for example, see (13.6.3) in \cite{AS}). Substituting $t=y^2$ in the integral, we obtain
\begin{align*}
&\int_{0}^\infty \calM_{s,0}(4\pi m y) \exp\left(-\frac{\pi N n^2 y^2}{v}\right)\frac{dy}{y}\\
&= 2^{2s-1}\Gamma(s+1/2)\int_{0}^\infty \sqrt{4\pi m y}I_{s-1/2}(2\pi m y)\exp\left(-\frac{\pi N n^2 y^2}{v}\right)\frac{dy}{y}\\
&= 2^{2s-1}\Gamma(s+1/2)\sqrt{\pi m}\int_{0}^\infty I_{s-1/2}(2\pi m \sqrt{t}) \exp\left(-\frac{\pi N n^2 t}{v}\right)t^{-3/4}\,dt.
\end{align*}
The latter integral is a Laplace transform which is computed in \cite{B2} (see (20) on p.197).
Inserting the evaluation, we obtain
\begin{align*}
&\int_{0}^\infty \calM_{s,0}(4\pi m y) \exp\left(-\frac{\pi N n^2 y^2}{v}\right)\frac{dy}{y}\\
&=2^{2s-1}\Gamma\left(s/2\right)\left(\frac{Nn^2}{\pi m^2 v}\right)^{1/4} M_{1/4,s/2-1/4}\left(\frac{\pi m^2 v}{Nn^2}\right)
\exp\left( \frac{\pi m^2 v}{2Nn^2}\right)\\
&=2^{2s-1}\Gamma\left(s/2\right)\left(\frac{Nn^2}{\pi m^2 v}\right)^{1/2}  \calM_{s/2+1/4 ,-1/2}\left(\frac{\pi m^2 v}{Nn^2}\right)
\exp\left( \frac{\pi m^2 v}{2Nn^2}\right).
\end{align*}
% Checked with Maple!!!
Consequently, we have in the case $n\mid m$ that
\begin{align*}
I(\tau,s,m,n)&=  \frac{2^{2s-2} Nn}{\sqrt{\pi}m}
\Gamma\left(s/2\right)
\calM_{s/2+1/4 ,-1/2}\left(\frac{\pi m^2 v}{Nn^2}\right)
 e\left(- \frac{m^2}{4Nn^2} u \right)\frake_{-m/n}.
\end{align*}
Substituting this in  \eqref{eq:lifti}, we find
\begin{align*}
\Lambda\big(\tau,F_{m}(z,s,-2)\big)&=\frac{2^{2-s}\sqrt{\pi} Ns(1-s)}{\Gamma(\frac{s}{2}-\frac{1}{2})}
\sum_{n\mid m} n\cdot
\calF_{\frac{m^2}{4Nn^2},-\frac{m}{n}}(\tau,\tfrac{s}{2}+\tfrac{1}{4},-\tfrac{1}{2}).
%\frac{1}{2\Gamma(2s)}\calM_{s/2+1/4 ,-1/2}\left(\frac{\pi m^2 v}{Nn^2}\right)
% e\left(- \frac{m^2}{4Nn^2} u \right)\frake_{-m/n}
%\mid_{-1/2,\rho_K}\gamma.
\end{align*}
Since $\calF_{m,h}(\tau,s,-1/2)=\calF_{m,-h}(\tau,s,-1/2)$, this concludes the proof of the theorem.
\end{proof}
%The following theorem is the main result of this section.

\begin{corollary}\label{cor:lift}
If $f\in H_{-2}(N)$ is a harmonic Maass form of weight $-2$ for $\Gamma_0(N)$,
then $\Lambda(\tau,f)$ belongs to $H_{-1/2,\rho_L}$. In particular, we have
\begin{align*}
\Lambda\big(\tau,F_{m}(z,2,-2)\big)=-2N
\sum_{n\mid m} n\cdot
\calF_{\frac{m^2}{4Nn^2},\frac{m}{n}}(\tau,\tfrac{5}{4},-\tfrac{1}{2}).
\end{align*}
\end{corollary}

\begin{proof}
The formula for the image of the Poincar\'e series $F_{m}(z,2,-2)$ is a direct consequence of Theorem \ref{thm:thetamain}. These Poincar\'e series for $m\in \Z_{>0}$ span the subspace $H^\infty_{-2}(N)\subset H_{-2}(N)$ of harmonic Maass forms whose principal parts at all cusps other than $\infty$ are constant.
Consequently, we find that the image of $H^\infty_{-2}(N)$ is contained   in $H_{-1/2,\rho_L}$.

For simplicity, here we only prove that the image of the full space $H_{-2}(N)$ is contained   in $H_{-1/2,\rho_L}$ in the special case when $N$ is squarefree. For general $N$ one can argue similarly, but the technical details get more complicated.
%
%Note that $\Orth(L'/L)$ can be identified with the group generated by the Atkin-Lehner involutions.
When $N$ is squarefree, then the group $\Orth(L'/L)$ of Atkin-Lehner involutions acts transitively on the cusps of $\Gamma_0(N)$. Consequently, we have
\[
H_{-2}(N)=\sum_{\gamma\in\Orth(L'/L)} \gamma H^\infty_{-2}(N).
\]
Using Proposition \ref{prop:acto},
%the equivariance of the theta  lift with respect to the action of $\Orth(L'/L)$,
we see that the whole space $H_{-2}(N)$ is mapped to $H_{-1/2,\rho_L}$.
\end{proof}

\begin{theorem}
The theta lift $\Lambda$ maps
%$M^!_{-2}(N)$ to $M^!_{-1/2,\rho_L}$.
weakly holomorphic modular forms to weakly holomorphic modular forms.
\end{theorem}

\begin{proof}
For simplicity we prove this only for the subspace $M^{!,\infty}_{-2}(N)$. If $N$ is squarefree, one obtains the result for the full space $M^!_{-2}(N)$ using the action of $\Orth(L'/L)$ as in the proof of Corollary \ref{cor:lift}. For general $N$, the argument gets more technical and we omit the details.

Let $F\in M^{!,\infty}_{-2}(N)$ and denote the Fourier expansion of $F$ at the cusp $\infty$ by
\[
F(z)=\sum_{m\in\Z} a_F(m)e(mz).
\]
We may write $F$ as a linear combination of Poincar\'e series as
\[
F(z)=\sum_{m>0}a_F(-m) F_m(z,2,-2).
\]
According to Corollary \ref{cor:lift}, we find that the principal part of $\Lambda(F)$ is equal to
\[
-2N\sum_{m>0}a_F(-m)
\sum_{n\mid m} n\cdot
e\left(-\frac{m^2}{4Nn^2}z\right)\left(\frake_{m/n}+\frake_{-m/n}\right).
\]
We now use the pairing $\{\cdot,\cdot\}$ of $H_{-1/2,\rho_L}$ with the space of cusp forms $S_{5/2,\bar\rho_L}$ (see \cite{BrF1}, Proposition 3.5) to prove that $\Lambda(F)$ is weakly holomorphic. We need to show that $\{\Lambda(F),g\}=0$ for every cusp form $g\in S_{5/2,\bar\rho_L}$. If we denote the coefficients of $g$ by $b(M,h)$, we
have
\begin{align*}
\{\Lambda(F),g\}&=-4N\sum_{m>0}a_F(-m)
\sum_{n\mid m} n\cdot b\left(\frac{m^2}{4Nn^2},\frac{m}{n}\right)\\
&= -4N\{F,\calS_1(g)\}.
\end{align*}
Here $\calS_1(g)\in S_4(N)$ denotes the (first) Shimura lift of $g$ as in \cite{Sk1}. Since $F$ is weakly holomorphic, the latter quantity vanishes.
\end{proof}

\begin{theorem}\label{GeneralThm}
Let  $f\in H_{-2}(N)$ and put $\partial f:= \frac{1}{4\pi} R_{-2,z}f$. For $m\in \Q_{>0}$ and $h\in L'/L$ the $(m,h)$-th Fourier coefficient of the
holomorphic part of $\Lambda(\tau,f)$ is equal to
\[
\tr_f(m,h)=-\frac{1}{2 m}\sum_{z\in Z(m,h)} \partial f(z).
\]
\end{theorem}

\begin{proof}
%\texttt{One can argue as in Katok-Sarnak \cite{KS}.}
%We put $\hat f:= \frac{1}{4\pi} R_{-2,z}f$.
%It follows from \eqref{eq:comm} that $\Delta_{0} \hat f=-2 \hat f$.
Inserting the definition of the theta lifting and using \eqref{eq:diffid}, we have
\begin{align*}
\Lambda\big(\tau,f\big)&=4\pi L_{3/2,\tau}\int_{M} \partial f(z) \Theta_L(\tau,z,\varphi_{KM})\\
&=\int_{M} \partial f(z) \Delta_{0,z}\Theta_L(\tau,z,\varphi_{S})\Omega.
\end{align*}
For $X\in V(\R)$ and $z\in \D$ we define $\varphi_S^0(X,z)=e^{2\pi Q(X)}\varphi_S(X,z)$. Then the Fourier expansion of the Siegel theta function in the variable $\tau$ is given by
\begin{align}
\label{eq:ft}
\Theta_L(\tau,z,\varphi_{S})= \sum_{X\in L'} \varphi_S^0(\sqrt{v}X,z)q^{Q(X)}\frake_X.
\end{align}

For $m\in \Q_{>0}$ and $h\in L'/L$, we put
$L_{m,h}=\{X\in L+h;\; Q(X)=m\}$. The group $\Gamma$ acts on $L_{m,h}$ with finitely many orbits.
We write $C(m,h)$ for the $(m,h)$-th Fourier coefficient of the
holomorphic part of $\Lambda(\tau,f)$. Using \eqref{eq:ft}, we see that
\[
C(m,h)=  \int_{M} \partial f(z) \Delta_{0,z} \sum_{X\in L_{m,h}}\varphi_{S}^0(\sqrt{v}X,z)\Omega.
\]
According to \cite[Proposition 3.2]{BrF1}, for $Q(X)>0$ the function $\varphi_{S}^0(X,z)$ has square exponential decay as $y\to \infty$. This implies that we may move the Laplacian in the integral to the function
$\partial f$. Since $\Delta_{0} \partial f=-2 \partial f$, we see that
\[
C(m,h)= -2\int_{M} \partial f(z) \sum_{X\in L_{m,h}}\varphi_{S}^0(\sqrt{v}X,z)\Omega.
\]
Using the usual unfolding argument, we obtain
\begin{align}
\label{eq:C}
C(m,h)= -2 \sum_{X\in \Gamma\bs L_{m,h}}\frac{1}{|\bar\Gamma_X|}\int_{\D} \partial f(z) \varphi_{S}^0(\sqrt{v}X,z)\Omega.
\end{align}

It is convenient to rewrite the integral over $\D$ as an integral over $G(\R)=\Sl_2(\R)$. If we normalize the Haar measure such that the maximal compact subgroup $\SO(2)$ has  volume $1$, we have
\begin{align*}
I(X):=\int_{\D} \partial f(z) \varphi_{S}^0(\sqrt{v}X,z)\,\Omega=
\int_{G(\R)} \partial f(gi) \varphi_{S}^0(\sqrt{v}X,gi)\,dg.
\end{align*}
Using the Cartan decomposition of $G(\R)$ and the uniqueness of spherical functions, we find, arguing as in the work of Katok and Sarnak \cite[pp.208]{KS}, that
\begin{align*}
I(X)=\partial f(\D_X)\cdot Y_\lambda(\sqrt{mv/N}),
\end{align*}
where
\[
Y_\lambda(t)=4\pi \int_1^\infty \varphi_{S}^0\left(t \alpha(a)^{-1} X(i),\,i\right) \omega_\lambda(\alpha(a))\,\frac{a^2-a^{-2}}{2}\,\frac{da}{a}.
\]
Here $\omega_\lambda(g)$ is the standard spherical function with eigenvalue $\lambda=-2$ (see e.g. \cite{Lang}, Chapters 5.4, 7.2, and 10.3), and $\alpha(a)=\kzxz{a}{0}{0}{a^{-1}}$.
%
% Note that the factor 4\pi is missing in [KS].
%
Note that $\omega_{-2}(\alpha(a))= \frac{a^2+a^{-2}}{2}$.
It is easily computed that
\[
\varphi_{S}^0\left(t \alpha(a)^{-1} X(i),\,i\right)
%= v e^{2\pi N t^2}e^{-\pi N t^2(a^4+a^{-4})},
=v e^{-\pi N t^2(a^2-a^{-2})^2},
\]
and therefore
\begin{align*}
Y_\lambda(t)&=2\pi v \int_0^\infty e^{-4\pi Nt^2\sinh(r)^2}\cosh(r)\sinh(r)\,dr
%&=2\pi v \int_0^\infty e^{-4\pi Nt^2 u^2}u\,du\\
=\frac{v}{4 Nt^2}.
\end{align*}
Hence
$Y_\lambda(\sqrt{mv/N})= \frac{1}{4 m}$.
%I(X&)=f'(\D_X).
Inserting this into \eqref{eq:C}, we obtain the assertion.
%\[
%C(m,h)= -\frac{1}{32\pi^2m} \sum_{X\in \Gamma\bs L_{m,h}}\frac{1}{|\bar\Gamma_X|} f'(\D_X).
%\]
\end{proof}

%\texttt{This seems to be off by a factor of $1/4\pi$.}

\begin{remark}
We can define similar theta liftings for other weights. For $k\in \Z_{\geq 0}$ odd we define a theta lifting of weak Maass forms of weight $-2k$ to weak Maass forms of weight $1/2-k$ by
\begin{align}
\label{eq:genlift1}
\Lambda(\tau,f,-2k)=  (L_{\tau})^{\frac{k+1}{2}}\int_{M} (R_{z}^{k} f)(z) \Theta_L(\tau,z,\varphi_{KM}).
\end{align}
In view of Proposition \ref{prop:eigval} and identity \eqref{eq:comm}, we should have that the lifting takes $H_{-2k}(N)$ to $H_{1/2-k,\rho_L}$ and maps weakly holomorphic forms to weakly holomorphic forms.
Analogously, for $k\in \Z_{\geq 0}$ even we define a theta lifting of weak Maass forms of weight $-2k$ to weak Maass forms of weight $3/2+k$ by
\begin{align}
\label{eq:genlift2}
\Lambda(\tau,f,-2k)=  (R_{\tau})^{\frac{k}{2}}\int_{M} (R_{z}^{k} f)(z)\Theta_L(\tau,z,\varphi_{KM}).
\end{align}
The lifting should then take $H_{-2k}(N)$ to $H_{3/2+k,\rho_L}$.
%and maps weakly holomorphic forms to weakly holomorphic forms.
These maps should give the interpretations in terms of a theta lift of the results discussed
in \cite[\S9]{Zagier}. Moreover, they should yield generalizations to congruence subgroups, arbitrary weights, and to harmonic Maass forms at the same time.
The lifting considered in \cite{BrF} is \eqref{eq:genlift2} in the case $k=0$. The lifting considered in the present paper is
%up to a constant factor
\eqref{eq:genlift1} for $k=1$.

\end{remark}

\subsection{The case of the partition function}
\label{sect:3.2}

Here  we  derive the formula for the partition function stated in Theorem \ref{main}  from Theorem \ref{GeneralThm} and Corollary \ref{cor:lift}.
We consider the theta lift of Section \ref{sect:3.1} in the special case when $N=6$.
We identify the discriminant form $L'/L$ with $\Z/12\Z$ together with the $\Q/\Z$-valued quadratic form $r\mapsto -r^2/24$.

The function $\eta(\tau)^{-1}$ can be viewed as a component of a vector valued modular form in $M^!_{-1/2,\rho_L}$ as follows. (Note that the latter space is isomorphic to the space $J_{0,6}^{weak}$ of weak Jacobi forms of weight $0$ and index $6$.)
We define
\[
G(\tau):=\sum_{r\in \Z/12\Z} \chi_{12}(r)\eta(\tau)^{-1}\frake_r.
\]
Using the transformation law of the eta-function under $\tau\mapsto \tau+1$ and $\tau\mapsto -1/\tau$, it is easily checked that
$G\in M^!_{-1/2,\rho_L}$. The principal part of $G$ is equal to $q^{-1/24}(\frake_{1}-\frake_{-5}-\frake_{7}+\frake_{11})$.

On the other hand, $G$  can be obtained as a theta lift. Let $F\in M^!_{-2}(6)$ be the function defined in \eqref{eq:defF}. It is invariant under the Fricke involution $W_6$, and under the Atkin-Lehner involution $W_3$ it is taken to its negative. Hence, in terms of Poincar\'e series we have
\[
F= F_1(\cdot ,2,-2) -F_1(\cdot ,2,-2)\mid W_2-F_1(\cdot ,2,-2)\mid W_3+ F_1(\cdot ,2,-2)\mid W_6.
\]
The function $P$ is given by $\frac{1}{4\pi}R_{-2}(F)$.
Using Corollary \ref{cor:lift} and Proposition \ref{prop:acto}, we see that $\Lambda(\tau,F)$ is an element of  $M^!_{-1/2,\rho_L}$ with principal part
$-4N q^{-1/24}(\frake_{1}-\frake_{-5}-\frake_{7}+\frake_{11})$. Consequently, we have
\[
G=-\frac{1}{4N}\cdot \Lambda(\tau,F).
\]
Now Theorem \ref{GeneralThm} tells us that for any positive integer $n$ the coefficient of $G$ with
index $(\tfrac{24n-1}{24},1)$ is equal to
\[
%tr_f(m,h)=
\frac{3}{N (24n-1)}\sum_{z\in Z\left(\tfrac{24n-1}{24},1\right)} P(z)
%= \frac{1}{ 24 n-1}\sum_{Q\in \calQ_{24n-1,1,6}/\Gamma_0(6)} P(\alpha_Q).
= \frac{1}{ 24 n-1}\sum_{Q\in \calQ_{n}} P(\alpha_Q).
\]
On the other hand, this coefficient is equal to $p(n)$ because
$$
\frac{q^{\frac{1}{24}}}{\eta(z)}=\prod_{n=1}^{\infty}\frac{1}{1-q^n}=
\sum_{n=0}^{\infty}p(n)q^n.
$$

%%%%%%%%%%%%%%%%%%%%%%%%%%%%%%%
\section{Complex Multiplication and singular moduli}\label{CM}
We have proved that
$p(n)=\Tr(n)/(24n-1)$.
To complete the proof of Theorem~\ref{main}, we require  results
from the theory of complex multiplication, and some new general results which
bound the denominators of singular moduli.

\subsection{Singular moduli for $j(z)$}
We first recall classical facts about Klein's $j$-function
\begin{equation}
j(z)=q^{-1}+744+196884q+21493760q^2+\dots.
\end{equation}
A point $\tau\in \H$ is a CM point if it is a root of a quadratic equation over
$\Z$. The {\it singular moduli} for $j(z)$, its values at such CM points, play a central role
in the theory of complex multiplication.
The following classical theorem (for example, see \cite{Borel, Cox}) summarizes
some of the most important properties of these numbers.

\begin{theorem}\label{CMtheorem}
Suppose that $Q=ax^2+bxy+cy^2$ is a primitive positive definite binary
quadratic form with discriminant $D=b^2-4ac<0$, and let $\alpha_Q\in \H$
be the point for which $Q(\alpha_Q,1)=0$. Then the following are true:
\begin{enumerate}
\item We have that $j(\alpha_Q)$ is an algebraic integer, and
 its minimal polynomial has degree $h(D)$,
the class number of discriminant $D$ positive definite binary quadratic forms.
\item The Galois orbit of $j(\alpha_Q)$ consists of the $j(z)$-singular moduli
associated to the $h(D)$ classes of discriminant $D$  forms.
\item If $K=\Q(\sqrt{D})$, then the discriminant $D$ singular moduli are conjugate to one another over $K$. Moreover, $K(j(\alpha_Q))$ is the discriminant $-D$
    Hilbert class field of $K$.
\end{enumerate}
\end{theorem}

Theorem~\ref{CMtheorem} and the properties of the weight 2 nonholomorphic
Eisenstein series
\begin{equation}\label{E2}
E_2^{*}(z):=-\frac{3}{\pi y}+E_2(z)=
1-\frac{3}{\pi y}-24\sum_{n=1}^{\infty}
\sum_{d\mid n}dq^n
\end{equation}
will play a central role in the proof of Theorem~\ref{main}.

\subsection{Bounding the denominators}\label{sect:bd}

Here we show that singular moduli like $6D\cdot P(\alpha_Q)$
are algebraic \emph{integers}, where $-D$ denotes the discriminant of $Q$.
We first introduce notation.
For a positive integer $N$, we let $\zeta_N$ denote a primitive $N$-th root of unity.
For a discriminant $-D<0$ and $r\in \Z$ with $r^2\equiv -D\pmod{4N}$ we let $\calQ_{D,r,N}$ denote the set of positive definite integral binary quadratic forms $[a,b,c]$ of discriminant $-D$ with $N\mid a$ and $b\equiv r\pmod{2N}$.
This notation is not to be confused with $\calQ_n$ introduced earlier. This set-up is more natural in this section.
For $Q=[a,b,c]\in \calQ_{D,r,N}$ we let $\alpha_Q=\frac{-b+\sqrt{-D}}{2a}$ be the corresponding {\em Heegner point} in $\H$.  We write $\calO_D$ for the order of discriminant $-D$ in $\Q(\sqrt{-D})$.

%The main result of this section is the following theorem.

\begin{theorem}
\label{thm:bound}
Let $D>0$ be coprime to $6$ and $r\in \Z$ with $r^2\equiv -D\pmod{24}$.
%Assume that $(D,6)=1$.
If $Q\in \calQ_{D,r,6}$ is primitive, then $6D\cdot P(\alpha_Q)$ is an algebraic integer contained in the ring class field corresponding to the order  $\calO_D\subset \Q(\sqrt{-D})$.
\end{theorem}

\begin{remark}  By Theorem~\ref{CMtheorem}, the multiset of  values $P(\alpha_Q)$
is a union of Galois orbits. Therefore, Theorem~\ref{thm:bound} completes the proof of Theorem~\ref{main}.
\end{remark}

Theorem \ref{thm:bound} will follow from Theorem \ref{thm:boundgen} below, a general result on values of derivatives of weakly holomorphic modular forms at Heegner points.
The following lemma is our key tool.

\begin{lemma}\label{bound2}
Let $\Gamma\subset \Gamma(1)$ be a level $N$ congruence subgroup.
Suppose that $f(z)$ is a weakly holomorphic modular function for $\Gamma$
whose Fourier expansions at all cusps have coefficients in $\Z[\zeta_N]$.
%\texttt{It seems that ``integral algebraic Fourier coefficients'' would be sufficient}.
If $\tau_0\in \H$ is a CM point, then $f(\tau_0)$ is an algebraic integer whose degree over
$\Q(\zeta_N,j(\tau_0))$ is bounded by $[\Gamma(1):\Gamma]$.
\end{lemma}

\begin{proof}
We consider the polynomial
\[
\Psi_f(X,z) = \prod_{\gamma\in \Gamma\bs \Gamma(1)} (X-f(\gamma z)).
\]
It is a monic polynomial in $X$ of degree $[\Gamma(1):\Gamma]$ whose coefficients are weakly holomorphic modular functions in $z$ for the group $\Gamma(1)$. Consequently, $\Psi_f(X,z)\in \C[j(z),X]$.

The assumption on the expansions of $f$ at all cusps means that for every $\gamma\in \Gamma(1)$ the modular function $f\mid \gamma$ has a Fourier expansion with coefficients in $\Z[\zeta_N]$.
So the coefficients of $\Psi_f(X,z)$ as a polynomial in $X$ are weakly holomorphic modular functions for $\Gamma(1)$ with coefficients in $\Z[\zeta_N]$, and therefore they are elements of $\Z[\zeta_N,j(z)]$.
Hence we actually have that $\Psi_f(X,z)\in \Z[\zeta_N,j(z),X]$.

Since $\Psi_f(f(z),z)=0$, we have, for every $z\in \H$, that $f(z)$ is integral over
$\Z[\zeta_N,j(z)]$ with degree bounded by $[\Gamma(1):\Gamma]$.
When $\tau_0$ is a CM point, then $j(\tau_0)$ is an algebraic integer, and the claim follows.
\end{proof}

\begin{comment}
\begin{remark}
If $f\in M_0^!(N)$ is defined over $\Q$ and $\tau_0$ is a Heegner point of level $N$, then it is well known that $f(\tau_0)\in \Q[j(\tau_0)]$. If in addition the Fourier expansion of $f$ at every cusp has coefficients in $\Z[\zeta_N]$, then Lemma \ref{bound2} implies that $f(\tau_0)$ is an algebraic integer
in $\Q[j(\tau_0)]$.
\end{remark}
\end{comment}

\subsection{Square-free level}

If the level $N$ is square-free, then the group of Atkin-Lehner involutions acts transitively on the cusps of $\Gamma_0(N)$. Combining this fact with Lemma \ref{bound2} leads to a handy criterion for the integrality of CM values. We begin by recalling some facts on Atkin-Lehner involutions
(see e.g. \cite{Kn} Chapter IX.7).

Let $N$ be an integer and $k$ a positive integer.  If $f$ is a complex valued function on the upper half plane $\H$ and $M=\kabcd\in \Gl_2^+(\R)$ then we put
\[
(f\mid M)(z) = (f\mid_k M)(z) = \det(M)^{k/2} (cz +d)^{-k} f(Mz).
\]
So scalar matrices act trivially.
%We extend the action to the group ring $\C[\Gl_2^+(\R)]$ in the usual way.
We write $M^!_k(N)$ for the space of weakly holomorphic modular forms of weight $k$ for the group $\Gamma_0(N)$.

\begin{comment}
The group $\Gamma_0(N)$ has index $N\prod_{p\mid N}(1+1/p)$ in $\Sl_2(\Z)$ and
$\sum_{0<d\mid N} \phi((d,N/d))$ equivalence classes of cusps.
The invariants of a cusp $a/c$ with $(a,c)=1$ are $d=(c,N)$ and $(c/(c,N))^{-1}a$ considered as element in the group of units of $\Z/(d,N/d)\Z$.
The width of the cusp $\kappa=a/c$ (where $(a,c)=1$) is $h=N/(c^2,N)$. Thus, if $M\in \Sl_2(\Z)$ with $M\infty=\kappa$, then $f\mid M$ has a Fourier expansion
\begin{align}\label{eq:fe}
f\mid M = \sum_{n\geq 0} a_\kappa(n) q^{n/h}.
\end{align}
\end{comment}

%We recall the Atkin-Lehner involutions on $M_k(N)$ (\cite{Kn} Chapter IX.7).
Let $Q$ be an exact divisor of $N$ (i.e. $Q|N$ and $(Q, N/Q)=1$), and let $W_Q^N$ be an integral matrix of the form
\[
W_Q^N=\zxz{Q\alpha}{\beta}{N\gamma}{Q\delta}
\]
with determinant $Q$.
If $f\in M^!_k(N)$, then $f\mapsto f|W_Q^N$ is independent of the choices of $\alpha,\beta,\gamma,\delta$, and defines an involution of $M^!_k(N)$, called an Atkin-Lehner involution.
If we write
\[
R_Q^N=\zxz{\alpha}{\beta}{N\gamma/Q}{Q\delta},
\]
we have $W_Q^N= R_Q^N \kzxz{Q}{0}{0}{1}$, and $R_Q^N\in\Gamma_0(N/Q)$.
For another   exact divisor $Q'$ of $N$, we have
\begin{align}
\label{eq:prod}
f\mid W_Q^N \mid W_{Q'}^N= f\mid W_{Q*Q'}^N,
\end{align}
where $Q*Q'=QQ'/(Q,Q')^2$.
If $(N',Q)=1$, then
\begin{align}
\label{eq:level}
f\mid W_Q^{NN'} =f \mid W_{Q}^{N}.
\end{align}
Clearly $W^N_N$ acts as the usual Fricke involution $W_N$.

From now on we assume that $N$ is square-free. Then the cusps of the group $\Gamma_0(N)$
are represented by $1/Q$, where $Q$ runs through the divisors $N$. Two cusps $a/c$ and $a'/c'$ (where $a,c,a',c'\in \Z$ and $(a,c)=(a',c')=1$) are equivalent under $\Gamma_0(N)$ if and only if $(c,N)=(c',N)$ (for example, see \cite{DS}, Prop.~3.8.3 and p.~103).
In particular, a complete set of representatives for the cusps of $\Gamma_0(N)$ is given by
$W_Q^N\infty$ with $Q$ running though the divisors of $N$.
Moreover, we have the disjoint  left coset decomposition
\begin{align}
\label{eq:decomp}
\Gamma(1) = \bigcup_{Q\mid N} \bigcup_{j\;(Q)} \Gamma_0(N) R_Q^N \zxz{1}{j}{0}{1}.
\end{align}

\begin{lemma}\label{bound3}
Let $N$ be square-free, and suppose that
$f\in M^!_0(N)$
%is a weakly holomorphic modular function on $\Gamma_0(N)$
has the property that $f\mid W_Q^N$ has  coefficients in $\Z$ for every
%Atkin-Lehner involution $W_Q^N$ with
$Q\mid N$.
%\texttt{It seems that ``integral algebraic Fourier coefficients'' would be sufficient}.
If $\tau_0$ is a level $N$ Heegner point of discriminant $-D$,
%in the imaginary quadratic field $\Q(\sqrt{D})$
then $f(\tau_0)$ is an algebraic integer in the ring class field for the order $\calO_D\subset \Q(\sqrt{-D})$.
\end{lemma}

\begin{proof}
The assumption on $f$ implies that $f\in \Q(j,j_N)$. Therefore, by the theory of complex multiplication (see Theorem \ref{CMtheorem}), $f(\tau_0)$ is contained in the claimed ring class field.
Since $N$ is square-free, the cusps of $\Gamma_0(N)$ are represented by $W_Q^N\infty$ with $Q\mid N$. Consequently, Lemma \ref{bound2} implies that $f(\tau_0)$ is an algebraic integer.
\end{proof}

\subsection{CM values of derivatives of weakly holomorphic modular forms}

The goal of this section is to prove the following theorem which easily implies Theorem \ref{thm:bound}.

\begin{theorem}\label{thm:boundgen}
Let $N$ be a square-free integer, and suppose that
$f\in M^!_{-2}(N)$
has the property that $f\mid_{-2} W_Q^N$ has  coefficients in $\Z$ for every
$Q\mid N$. Define $\partial f=\frac{1}{4\pi} R_{-2} f$.
%\texttt{It seems that ``integral algebraic Fourier coefficients'' would be sufficient}.
%If $\tau_0$ is a level $N$ Heegner point of discriminant $-D$,
%in the imaginary quadratic field $\Q(\sqrt{D})$
%then $6D\cdot f'(\tau_0)$ is an algebraic integer in the ring class field for the order %$\calO_D\subset \Q(\sqrt{-D})$.
%
Let $D>0$ be coprime to $2N$ and $r\in \Z$ with $r^2\equiv -D\pmod{4N}$.
If $Q\in \calQ_{D,r,N}$ is primitive, then $6D\cdot \partial f(\alpha_Q)$
is an algebraic integer in the ring class field for the order  $\calO_D\subset \Q(\sqrt{-D})$.
\end{theorem}

To prove the theorem we need two lemmas.

\begin{lemma}
\label{lem:bound1new}
Let $N$ be a square-free integer.
Let $f\in M_{-2}^!(N)$ and assume that $f\mid_{-2}  W_Q^N$ has integral Fourier coefficients for all $Q\mid N$. Define
\[
A_f  = \frac{1}{4\pi} R_{-2} f-\frac{1}{6}f E_2^*.
\]
Let $D>0$ and $r\in \Z$ with $r^2\equiv -D\pmod{4N}$.
If $Q\in \calQ_{D,r,N}$ is primitive, then $6\cdot A_f(\alpha_Q)$
is an algebraic integer in the ring class field for the order  $\calO_D\subset \Q(\sqrt{-D})$.
\end{lemma}

\begin{proof}
First, computing the Fourier expansion, we notice that $A_f\in M_0^!(N)$.
Then, using the fact that $R_{-2}(f\mid_{-2} W^N_Q)= (R_{-2}f)\mid_{-2} W^N_Q$, we see that $6A_f\mid W_Q^N$ has integral coefficients   for all $Q\mid N$.
Consequently, the assertion follows from Lemma \ref{bound3}.
\end{proof}

\begin{lemma}
\label{lem:bound2}
Let $N$ be a  square-free integer.
Let $f\in M_{-2}^!(N)$ and assume that $f\mid W_Q^N$ has integral Fourier coefficients for all $Q\mid N$. Define
$\hat f  = f\cdot E_2^*$.
Let $D>0$ be coprime to $2N$ and $r\in \Z$ with $r^2\equiv -D\pmod{4N}$.
If $Q\in \calQ_{D,r,N}$ is primitive, then $D\cdot \hat f(\alpha_Q)$
is an algebraic integer in the ring class field for the order  $\calO_D\subset \Q(\sqrt{-D})$.
\end{lemma}

\begin{proof}
We write $Q=[a,b,c]$. Since $-D=b^2-4ac$ is odd, $b$ is odd. Hence
\[
M=\zxz{-b}{-2c}{2a}{b}
\]
is a primitive integral matrix of determinant $D$, satisfying $M\alpha_Q= \alpha_Q$.
By the elementary divisor theorem there exist $\gamma_1, \gamma_2\in \Gamma_0(N)$ such that
\begin{align}
\label{eq:eldiv}
M=\gamma_1^{-1} \zxz{1}{0}{0}{D}\gamma_2.
\end{align}
% See e.g. Miyake Lemma 4.5.2.
We put
\[
r_D(z) = E_2^*(z)-DE_2^*(Dz) = E_2^*(z)-(E_2^*\mid W_D)(z).
\]
Because of \eqref{eq:eldiv}, we have
\[
E_2^*\mid M = r_D\mid \gamma_1 M + E_2^*.
\]
Using the fact that $(2a\alpha_Q+b)^2=-D$, we find that
%%arguing as in \cite{millerpixton}, Section 3.2,
 that
\begin{align*}
E_2^*(\alpha_Q) &= \frac{1}{2} (r_D\mid_2 \gamma_1)(\alpha_Q),\\
\hat f (\alpha_Q) & = \frac{1}{2} (f \cdot r_D)(\gamma_1\alpha_Q).
\end{align*}
Arguing as in the proof of Proposition 3.1 of \cite{millerpixton}, we see that $\hat f(\alpha_Q)$ is contained in the claimed ring class field. Hence,
replacing $Q$ by $\gamma_1 Q\in \calQ_{D,r,N}$, it suffices to prove
that $\frac{D}{2} (f \cdot r_D)(\alpha_Q)$ is an algebraic integer.
% in the claimed ring class field.

In view of Lemma \ref{bound2} it suffices to show that for any $\gamma\in \Gamma(1)$, the weakly holomorphic modular form
$\frac{D}{2}(f\cdot r_D)\mid \gamma$ has Fourier coefficients in $\Z[\zeta_{ND}]$.
According to \eqref{eq:decomp}, there exists a $\gamma'\in \Gamma_0(N)$, a divisor $Q\mid N$ and $j\in \Z$ such that
\[
\gamma = \gamma' R_Q^N \zxz{1}{j}{0}{1}= \gamma' W_Q^N \zxz{1/Q}{j/Q}{0}{1}.
\]
Consequently, we have
\[
(f\mid_{-2} \gamma)(z) = Q\cdot \left(f\mid_{-2} W_Q^N\right)\left( \frac{z+j}{Q}\right) \in Q\cdot\Z[\zeta_Q]((q^{1/Q})).
\]

To analyze the situation for $r_D$, we write the integral matrix $W_D \gamma$ of determinant $D$ as
\[
W_D \gamma =\gamma'' \zxz{D_1}{k}{0}{D_2}
\]
with $\gamma''\in \Gamma(1)$ and positive integers $D_1,D_2,k$ satisfying  $D_1 D_2=D$.
Then we have
\begin{align*}
(r_D\mid_2 \gamma)(z) &= (E_2^*\mid \gamma)(z) - (E_2^*\mid W_D \gamma)(z)\\
&= E_2^*(z) - \frac{D}{D_2^2} E_2^*\left( \frac{D_1 z+k}{D_2}\right) .
\end{align*}
Taking into account that $D$ is odd, we see that $\frac{D}{2}(r_D\mid_2 \gamma)\in \Z[\zeta_D]((q^{1/D}))$. This concludes the proof of the lemma.
\end{proof}

\begin{proof}[Proof of Theorem~\ref{thm:boundgen}]
Using the notation of Lemma \ref{lem:bound1new} and Lemma \ref{lem:bound2}, we have
\[
\partial f(z)=A_f(z)+\frac{1}{6}\hat f(z).
\]
Consequently, the assertion follows from these lemmas.
\end{proof}

\begin{proof}[Proof of Theorem~\ref{thm:bound}]
We apply Theorem~\ref{thm:boundgen} to the function $F\in M_{-2}^!(6)$ defined in \eqref{eq:defF}.
Note that $F\mid W_6^6=F$ and $F\mid W_3^6= -F$. Moreover, we have $P=\frac{1}{4\pi}R_{-2} F=\partial F$.
\end{proof}
%%%%%%%%%%

%%%%%%%%%%%%%%%%%%%%%%%%%%%%%%
\section{Examples}\label{examples}

To compute $p(n)$ using Theorem~\ref{main}, one first determines
representatives for $\calQ_n=Q_{24n-1,1,6}$, a set which has
$h(-24n+1)$ many elements. Gross, Kohnen, and Zagier (see pages 504--504 of \cite{GKZ})
establish a one to one correspondence between representatives of
$\calQ_n$ and positive definite binary quadratic forms under $\SL_2(\Z)$
with discriminant $-24n+1$. Therefore, to determine representatives for $\calQ_n$,
it suffices to use the theory of {\it reduced forms} (for example, see page 29
of \cite{Cox})
to determine representatives for the $\SL_2(\Z)$ equivalence classes, and to then
apply the Gross-Kohnen-Zagier correspondence (see the Proposition on page 505 of
\cite{GKZ}).

For example, if $n=2$, then $-24n+1=-47$, and we have that $h(-47)=5$ and
$\calQ_2=\{[6,1,2],\ [12,1,1],\ [18,13,3],\ [24,25,7],\ [36,49,17]\}$,
where $[a,b,c]:=ax^2+bxy+cz^2$.
Calculating $p(n)$, by Theorem~\ref{main}, now follows
from sufficiently accurate numerical approximations of the algebraic integers $6(24n-1)P(\alpha_Q)$.

We used this method to compute
the first few ``partition polynomials'' $H_n(x)$.

\bigskip
\begin{center}
\begin{tabular}{cccccc}
\hline\hline
$n$ && $(24n-1)p(n)$ && $H_n(x)$
\\ \hline \\
$1$ &&  $23$ && $x^3-23x^2+\frac{3592}{23}x-419$  \\ \ \ \\
$2$ &&  $94$ && $x^5-94x^4+\frac{169659}{47}x^3-65838x^2+\frac{1092873176}{47^2}x
+\frac{1454023}{47}$  \\ \ \ \\
$3$ && $213$ && $x^7-213x^6+\frac{1312544}{71}x-723721x^4+\frac{44648582886}{71^2}x^3$\\
\ \ \\
\ \ && \ \ &&\ \ \ \ \ \ \ \ \ \ \ \ \ $+\frac{9188934683}{71}x^2+
\frac{166629520876208}{71^3}x+\frac{2791651635293}{71^2}$\\
\ \ \\
$4$ && $475$ && $x^8-475x^7+\frac{9032603}{95}x^6-9455070x^5+
\frac{3949512899743}{95^2}x^4$\\ \ \ \\
\ \ && \ \ && \ \ \ \ \ \ \ $-\frac{97215753021}{19}x^3+\frac{9776785708507683}{95^3}x^2$
\\ \ \ \\
\ \ && \ \ && \ \ \ \ \ \ \ \ \ \ \ $-\frac{53144327916296}{19^2}x-
\frac{134884469547631}{5^4\cdot 19}$.\\
\end{tabular}
\end{center}

\bigskip

We conclude with some natural questions which merit further investigation.
\medskip

\begin{enumerate}
\item Is it true that $(24n-1)P(\alpha_Q)$ is an algebraic integer for all $Q\in \calQ_n$?
\item What can be said about the irreducibility of the $H_n(x)$?
\item Is there a ``closed formula'' for the constant terms of the $H_n(x)$
which is analogous to the formula of Gross and Zagier \cite{GrossZagier} on norms of differences of $j(z)$-singular moduli?
\item Do the singular moduli $P(\alpha_Q)$ enjoy special congruence properties? If so,
do such congruences imply Ramanujan's congruences modulo 5, 7, and 11?
\end{enumerate}


\begin{thebibliography}{GKZ}

\bibitem{AS} M. Abramowitz and I. Stegun, \emph{Pocketbook of Mathematical Functions}, Verlag Harri Deutsch (1984).

\bibitem{Andrews} G. E. Andrews, \emph{The theory of partitions}, Cambridge Univ. Press,
Cambridge, 1984.

\bibitem{Bo1}
R. Borcherds,
\emph{Automorphic forms with singularities on Grassmannians},
Invent. Math. \textbf{132} (1998), 491-562.

\bibitem{Borel} A. Borel, S. Chowla, C. S. Herz, K. Iwasawa, and J.-P. Serre,
\emph{Seminar on complex multiplication}, Springer Lect. Notes. \textbf{21},
Springer Verlag, Berlin, 1966.

%%\bibitem{AL} A. O. L. Atkin and J. Lehner,
%%\emph{Hecke operators on $\Gamma_0(m)$},
%%Math. Ann. \textbf{185} (1970),  134-160.

\bibitem{BringmannOno} K. Bringmann and K. Ono,
\emph{Arithmetic properties of coefficients of half-integral weight Maass-Poincar\'e series},
Math. Ann. \textbf{337} (2007),  591-612.

\bibitem{BringmannOnoPAMS} K. Bringmann and K. Ono, \emph{An arithmetic formula for the partition function},
Proc. Amer. Math. Soc. \textbf{135} (2007),  3507-3514.

\bibitem{Br} J. Bruinier, \emph{Borcherds products on
  $\Orth(2,l)$ and Chern classes of Heegner divisors},  Springer Lecture Notes in Mathematics {\bf 1780}, Springer-Verlag (2002).

\bibitem{BrF1} J. H. Bruinier and J. Funke,
\emph{On two geometric theta lifts},
Duke Math. J. \textbf{125} (2004),  45-90.

\bibitem{BrF} J. H. Bruinier and J. Funke,
\emph{Traces of CM-values of modular functions},
J. Reine Angew. Math. \textbf{594} (2006),
1-33.

\bibitem{BrJO} J. H. Bruinier, P. Jenkins, and K. Ono,
\emph{Hilbert class polynomials and traces of singular moduli},
Math. Ann. \textbf{334} (2006),  373-393.

\bibitem{Cox} D. A. Cox, \emph{Primes of the form $x^2+ny^2$},
Wiley and Sons, New York, 1989.

\bibitem{DS} F. Diamond and J. Shurman,  \emph{A first course in modular forms}, Graduate Texts in Mathematics {\bf 228}, Springer-Verlag (2005).

\bibitem{DIT1}  W. Duke, \"O. Imamo$\overline{{\text {\rm g}}}$lu, and A. Toth, \emph{Cycle integrals of the $j$-function and mock modular forms},
Ann. of Math., accepted for publication.

\bibitem{DIT2} W. Duke, \"O. Imamo$\overline{{\text {\rm g}}}$lu, and A. Toth, \emph{Real quadratic analogues of traces of singular invariants},
Int. Math. Res. Notices, accepted for publication.

\bibitem{DukeJenkins} W. Duke and P. Jenkins,  \emph{Integral traces of singular
values of weak Maass forms}, Algebra and Number Th. \textbf{2} (2008),
573-593.

\bibitem{B2} A. Erd\'elyi, W. Magnus, F. Oberhettinger and F. G. Tricomi, \emph{Tables of Integral Transforms}, vol.~I, McGraw-Hill (1954).

\bibitem{FolsomMasri} A. Folsom and R. Masri, \emph{Equidistribution of Heegner points
and the partition function}, Math. Ann., accepted for publication.

\bibitem{Funke} J. Funke, \emph{Heegner divisors and nonholomorphic modular forms},
Compositio Math. \textbf{133} (2002),  289-321.


\bibitem{GKZ} B. Gross, W. Kohnen, and D. Zagier, \emph{Heegner
  points and derivatives of $L$-series. II}.  Math. Ann.  {\bf 278}
  (1987), 497--562.

\bibitem{GrossZagier}  B. Gross and D. Zagier, \emph{On singular moduli}, J. Reine Angew. Math.
\textbf{355} (1985), 191-220.

\bibitem{HardyRamanujan} G. H. Hardy and S. Ramanujan, \emph{Asymptotic formulae in combinatory
analysis}, Proc. London Math. Soc. (2) \textbf{17} (1918), 75-115.

\bibitem{Jenkins} P. Jenkins, \emph{Kloosterman sums and traces of
singular moduli}, J. Number Th. \textbf{117} (2006), 301-314.

\bibitem{KS} S. Katok and P. Sarnak, \emph{Heegner points, cycles,
and Maass forms}, Israel J. Math. \textbf{84} (1993),
193-227.

\bibitem{Kn}  A. W. Knapp, \emph{Elliptic curves}, Princeton University Press (1992).

\bibitem{KAnn}
S. Kudla,
\emph{Central derivatives of Eisenstein series and height pairings},
Ann. of Math. \textbf{146} (1997), 545-646.

\bibitem{KMI}
S. Kudla and J. Millson,
\emph{The Theta Correspondence and Harmonic Forms I},
Math. Ann. \textbf{274} (1986),  353-378.

\bibitem{KM90}
S. Kudla and J. Millson,
\emph{Intersection numbers of cycles on locally symmetric spaces and Fourier coefficients of holomorphic modular forms in several complex variables},
IHES Pub. \textbf{71} (1990), 121-172.

\bibitem{Lang} S. Lang, \emph{$SL_2(\R)$}, Graduate Texts in Mathematics {\bf 105}, Springer, New York (1985).

\bibitem{millerpixton} A. Miller and A. Pixton,
\emph{Arithmetic traces of non-holomorphic modular invariants},
Int. J. Number Th. \textbf{6} (2010),  69-87.

\bibitem{cbms} K. Ono, \emph{The web of modularity: Arithmetic
of the coefficients of modular forms and $q$-series}, CBMS
Regional Conference, {\bf 102}, Amer. Math. Soc., Providence, R.
I., 2004.

\bibitem{Rademacher1} H. Rademacher, \emph{On the partition function $p(n)$},
Proc. London Math. Soc. (2) \textbf{43} (1937),  241-254.

\bibitem{Rademacher2} H. Rademacher, \emph{On the expansion of the partition function
in a series}, Ann. Math. \textbf{44} (1943),  416-422.

\bibitem{Sk1} N.-P. Skoruppa, \emph{Developments in the theory
  of Jacobi forms}, In: Proceedings of the conference on automorphic
  funtions and their applications, Chabarovsk (eds.: N. Kuznetsov and
  V. Bykovsky), The USSR Academy of Science (1990), 167--185. (see
  also MPI-preprint 89-40, Bonn (1989).)

\bibitem{Zagier} D. Zagier, \emph{Traces of singular moduli},
Motives, polylogarithms and Hodge theory, Part I (Irvine, CA, 1998)
(2002), Int. Press Lect. Ser., 3, I, Int. Press, Somerville, MA,
211--244.



\end{thebibliography}
\end{document}